\documentclass[reqno]{amsart}

\usepackage{amsthm, amssymb, latexsym, tabls, hhline,lscape}
\newtheorem{theorem}{Theorem}[section]

\newtheorem{lemma}[theorem]{Lemma}
\theoremstyle{definition}

\theoremstyle{remark}
\newtheorem{remark}[theorem]{Remark}
\numberwithin{equation}{section}
\newcommand{\Sym}{\ensuremath{\mathit{Sym}}}
\newcommand{\QSym}{\ensuremath{\mathit{QSym}}}
\newcommand{\NSym}{\ensuremath{\mathit{NSym}}}
\newcommand{\ra}{\rightarrow}
\newcommand{\ot}{\ensuremath{\otimes}}
\newcommand{\scal}[1]{\langle #1 \rangle}

\newcommand{\Z}{\mathbb Z}

\makeindex

\title[Combinatorial Hopf algebras and Towers of Algebras]
{\bf Combinatorial Hopf algebras\\ and Towers of Algebras}

\author{Nantel Bergeron}\address[Nantel Bergeron]
{Department of Mathematics and Statistics\\ York  University\\ To\-ron\-to, Ontario M3J 1P3\\ CANADA}
\email{bergeron@mathstat.yorku.ca}
\urladdr{http://www.math.yorku.ca/bergeron}

 \author{Thomas Lam}\address[Thomas Lam]
 {Department of Mathematics\\ Harvard University\\ Cambridge,\\ MA 02138.}
 \email{tfylam@math.harvard.edu}
 \urladdr{http://www.math.harvard.edu/\~{ }tfylam}

 \author{Huilan Li}\address[Huilan Li]
 {Department of Mathematics and Statistics\\ York  University\\ To\-ron\-to, Ontario M3J 1P3\\ CANADA}
 \email{lihuilan@mathstat.yorku.ca}

\date{\today}

\thanks{N.B. anf H. L. are supported in part by CRC and NSERC.}
\thanks{T.L. is partially supported by NSF grants DMS-0600677 and DMS-0652641.}



\begin{document}
\maketitle

%%%%%%%%%%%%%%%%%%%%%%%%%%%%%%%%%%%
%%%%%%%%%%%%%%%%%%%%%%%%%%%%%%%%%%%
%\begin{abstract}
%The Grothendieck group of the tower of symmetric group algebras
%has a graded self-dual Hopf algebra structure. Inspired by this Bergeron and Li introduce by means of axioms, a general notion of a tower of
%algebras $\bigoplus_{n\ge0}A_n$ and study two Grothendieck groups on this tower linked by
%a natural paring. They show that the
%axioms give a structure of graded Hopf algebras on each
%Grothendieck groups and these structures are graded dual to each other.
%In this paper we show the following very surprising fact.
%For any tower $\bigoplus_{n\ge0}A_n$, we must have that $\dim(A_n)=r^nn!$.
%\end{abstract}

\begin{abstract}
Bergeron and Li have introduced a set of axioms which guarantee that
the Grothendieck groups of a tower of
algebras $\bigoplus_{n\ge0}A_n$ can be endowed with the structure of graded dual Hopf
algebras.  Hivert and Nzeutzhap, and independently Lam and Shimozono
constructed dual graded graphs from primitive elements in Hopf
algebras.  In this paper we apply the composition of these
constructions to towers of algebras. We show that if a tower
$\bigoplus_{n\ge0}A_n$ gives rise to graded dual Hopf algebras then
we must have $\dim(A_n)=r^nn!$ where $r = \dim(A_1)$.
\end{abstract}

%%%%%%%%%%%%%%%%%%%%%%%%%%%%%%%%%%%
%%%%%%%%%%%%%%%%%%%%%%%%%%%%%%%%%%%
\section{Introduction}
\setcounter{equation}{0}

This paper is concerned with the interplay between towers of
associative algebras, pairs of dual combinatorial Hopf algebras, and
dual graded graphs.  Our point of departure is the study of the
composition of two constructions: (1) the construction of dual Hopf
algebras from towers of algebras satisfying some axioms, due to
Bergeron and Li \cite{BL}; and (2) the construction of dual graded
graphs from primitive elements in dual Hopf algebras, discovered
independently by Hivert and Nzeutchap \cite{HN}, and Lam and
Shimozono:
\begin{equation}\label{eq:comp}
\text{tower of algebras} \longrightarrow \text{combinatorial Hopf
algebra} \longrightarrow \text{dual graded graph}
\end{equation}
The notion of a pair $(\Gamma,\Gamma')$ dual graded graphs is likely
to be the least familiar.  They were introduced by Fomin \cite{F}
(see also \cite{S}) to encode the enumerative properties of the
Robinson-Schensted correspondence and its generalizations.  The
first arrow in (\ref{eq:comp}) is obtained by using induction and
restriction on the Grothendieck groups.  The second arrow is
obtained by using (some of the) structure constants of a
combinatorial Hopf algebra as edge multiplicities for a graph.  We
review these constructions in Sections \ref{sec:BL} and
\ref{sec:HNLS}. The notion of combinatorial Hopf algebra used here  is related, but slightly different from the one in  \cite{ABS}.

The key example of all three classes of objects arises from the
theory of symmetric functions.  In ~\cite{LG}, L. Geissinger showed
that the ring \Sym\ of symmetric functions is a graded self-dual
Hopf algebra.  Using the work of Frobenius and Schur,
Zelevinsky~\cite{AVZ} interpreted the Hopf structure in terms of the
Grothendieck groups of the tower of symmetric group algebras
$\bigoplus_{n\geq0}\mathbb{C}\mathfrak{S}_{n}$.  Finally it follows
from the classical work of Young that the branching rule for the
symmetric group, or equivalently the Pieri rule for symmetric
functions, gives rise to the {\it Young graph} on the set of
partitions.  The Young graph is the motivating example of dual
graded graphs.

In recent years it has been shown that other graded dual Hopf
algebras can be obtained from towers of algebras.  In~\cite{MR}
Malvenuto and Reutenauer establish the duality between the Hopf
algebra \NSym\ of noncommutative symmetric functions and the Hopf
algebra \QSym\ of quasi-symmetric functions. Krob and
Thibon~\cite{KT} then showed that this duality can be interpreted as
the duality of the Grothendieck groups associated with
$\bigoplus_{n\geq0}H_{n}(0)$ the tower of Hecke algebras at $q = 0$.
For more examples, see~\cite{BHT,HNT,ANS}.

It is very tempting, as suggested by J. Y. Thibon, to classify all
combinatorial Hopf algebras which arise as Grothendieck groups
associated with a tower of algebras $\bigoplus_{n\geq0}A_{n}$. The
list of axioms given by the first and last author in \cite{BL}
guarantees that the Grothendieck groups of a tower of algebras form
a pair of graded dual Hopf algebras.  This list of axioms is not
totally satisfactory as some of the axioms are difficult to verify
and the description is far from a classification. In this paper we
present a very surprising fact which shows that towers of algebras
giving rise to combinatorial Hopf algebras are much more rigid than
they appear.

\begin{theorem} \label{dimthm}
If $A=\bigoplus_{n\geq0}A_{n}$ is  a tower of algebras such that its
associated Grothendieck groups form a pair of graded dual Hopf
algebras, then $\dim(A_n)=r^nn!$ where $r=\dim(A_1)$.
\end{theorem}

The notion of ``forming a pair of graded dual Hopf algebras'' is
made precise in Section \ref{sec:BL}.  The numbers $r^n n!$ will be
familiar to experts in the theory of dual graded graphs -- they
count certain paths in a pair of dual graded graphs.

The rigidity proved in Theorem \ref{dimthm} suggests that there may
be a structure theorem for towers of algebras which give rise to
combinatorial Hopf algebras.  In particular, to perform the inverse
constructions of the arrows in (\ref{eq:comp}), it suggests that one
should study algebras related to symmetric groups (or wreath
products of symmetric groups).  There are many combinatorial Hopf
algebras for which one may attempt to perform the inverse
construction, but there are even more dual graded graphs.  The general
construction of \cite{LamSKac} produces dual graded graphs from
Bruhat orders of Weyl groups of Kac-Moody algebras and it is unclear
whether there are Hopf algebras, or towers of algebras giving rise
to these graphs.

{\bf Acknowledgements.}  The second author would like to thank Mark
Shimozono for the collaboration which led to the line of thinking in
this paper.

%Before we prove this theorem, we recall the definition of  a tower
%of algebras and its associated Grothendieck groups in Section~2.
%Once this is done, we prove the theorem in Section~3 using the
%theory of dual graded graphs of Fomin~\cite{F} and Stanley~\cite{S}.

%%%%%%%%%%%%%%%%%%%%%%%%%%%%%%%%%%%
%%%%%%%%%%%%%%%%%%%%%%%%%%%%%%%%%%%
\section{From towers of algebras to combinatorial Hopf algebras}
\label{sec:BL} We recall here the work of Bergeron and Li~\cite{BL}
on towers of algebras.  %For more details on the the theory of
%bialgebras, representations and Grothendieck groups
%see~\cite{ARS,CR,LG}.
For $B$ an arbitrary algebra we denote by $
_{B}\mbox{mod}$, the category of all finitely generated left
$B$-modules, and by $\mathcal{P}(B)$, the category of all finitely
generated projective left $B$-modules. For some category
$\mathcal{C}$ of left $B$-modules ($_{B}\mbox{mod}$ or
$\mathcal{P}(B)$) let $\mathbf{F}$ be the free abelian group
generated by the symbols $(M)$, one for each isomorphism class of
modules $M$ in $\mathcal{C}$. Let $\mathbf{F_{0}}$ be the subgroup
of $\mathbf{F}$ generated by all expressions $(M)-(L)-(N)$ one for
each exact sequence
              $$0\rightarrow L\rightarrow M\rightarrow N\rightarrow0$$
in $\mathcal{C}$. The
\textit{Grothendieck group} ${\mathcal K}_{0}(\mathcal{C})$ of the category
$\mathcal{C}$ is defined by the quotient $\mathbf{F}/\mathbf{F_{0}}$, an abelian additive group.
For $M\in\mathcal{C}$, we denote by $[M]$ its image in ${\mathcal K}_{0}(\mathcal{C})$.
We then set
             $$G_{0}(B)={\mathcal K}_{0}(_{B}\mbox{mod}) \quad\hbox{and}\quad
            K_{0}(B)={\mathcal K}_{0}(\mathcal{P}(B)).$$

For $B$ a finite-dimensional algebra over a field $K$, let
$\{V_{1},\cdots, V_{s}\}$ be a complete list of nonisomorphic
simple $B$-modules. The projective covers $\{P_{1},\cdots,
P_{s}\}$ of the simple modules $V_i$'s  is a complete list of nonisomorphic indecomposable
projective $B$-modules. We have that  $G_{0}(B)=\bigoplus^s_{i=1}\mathbb{Z}[V_{i}]$ and
$ K_{0}(B)=\bigoplus^s_{i=1}\mathbb{Z}[P_{i}]$.

Let $\varphi\colon B\to A$ be an injection of algebras preserving
unities, and let $M$ be a (left) $A$-module and $N$ a (left)
$B$-module. The \textit{induction} of $N$ from $B$ to $A$ is
Ind$^{A}_{B}N=A\otimes_{\varphi}N$, the (left) $A$-module $A\otimes
N$ modulo the relations $a\otimes bn\equiv a\varphi(b)\otimes n$,
and the \textit{restriction} of $M$ from $A$ to $B$ is
Res$^{A}_{B}M=\mbox{Hom}_{A}(A,M)$, the (left) $B$-module with the
$B$-action
 defined by $bf(a)=f(a\varphi(b))$.

Let $A=\bigoplus_{n\geq 0}A_{n}$ be a graded algebra over
$\mathbb{C}$ with multiplication $\rho\colon A\otimes A\to A$.
 Bergeron and Li studied five axioms for $A$ (we refer to \cite{BL} for full details):

\smallskip
\noindent(1) For each $n\ge 0$, $A_{n}$ is a finite-dimensional
algebra by itself with  (internal) multiplication $\mu_n\colon
A_n\otimes A_n\to A_n$  and unit $1_n$. $A_{0}\cong \mathbb{C}$.

\smallskip
\noindent(2) The (external) multiplication $\rho_{m,n}:
A_{m}\otimes A_{n}\rightarrow A_{m+n}$ is an injective
homomorphism of algebras, for all $m\mbox{ and }n$ (sending
$1_{m}\otimes 1_{n}$ to $1_{m+n}$ ).

\smallskip
\noindent(3) $A_{m+n}$ is a two-sided projective $A_{m}\otimes
A_{n}$-module with the action defined by $a\cdot(b\ot
c)=a\rho_{m,n}(b\ot c)\mbox{ and }(b\ot c)\cdot a=\rho_{m,n}(b\ot
c)a$, for all $m,n\geq0,\ a\in A_{m+n},\ b\in A_{m},\ c\in
A_{n}\mbox{ and }m,n\geq0$.

\smallskip
\noindent(4) A relation between the decomposition of $A_{n+m}$ as a
left $A_m\ot A_n$-module and as a right $A_m\ot A_n$-module holds.

\smallskip
\noindent(5) An analogue of Mackey's formula relating induction and restriction of modules holds.

\medskip
We say here that $A = \bigoplus_{n\geq 0}A_{n}$ is a {\it tower of
algebras} if it satisfies Conditions (1), (2) and (3).

Condition (1) guarantees that the Grothendieck groups
$G(A)=\bigoplus_{n\ge 0} G_0(A_n)$ and $K(A)=\bigoplus_{n\ge 0}
K_0(A_n)$ are graded connected.  Conditions (2) and (3) ensure that
induction and restriction are well defined on $G(A)$ and $K(A)$,
defining a multiplication and comultiplication, as follows.  For
$[M]\in G_0(A_m)$ (or $K_0(A_m)$) and $[N]\in G_0(A_n)$ (or
$K_0(A_n)$) we let
 $$[M] [N]=\left[\mbox{Ind}^{A_{m+n}}_{A_{m}\otimes
A_{n}}M\otimes N\right]  \qquad\hbox{and}\qquad
\Delta([N])=\sum_{k+l=n} \left[\mbox{Res}^{A_{k+l}}_{A_{k}\otimes
A_{l}} N\right].$$ The pairing between $K(A)$ and $G(A)$ is given by
$\scal{\,\,,\,}: K(A)\times G(A)\ra \mathbb{\mathbb{Z}} $
where
        $$\scal{[P],[M]}=\left\{\begin{array}{ll}
                       \mbox{dim}_{K}\big(\mbox{Hom}_{A_{n}}(P,M)\big)
                       & \mbox{if }[P]\in K_{0}(A_{n})
                       \mbox{ and }[M]\in G_{0}(A_{n}),\\
                       0 & \mbox{otherwise.}
                       \end{array}
                \right.$$

Thus with (only) Conditions (1), (2), and (3), $G(A)$ and $K(A)$ are
dual free $\Z$-modules both endowed with a multiplication and
comultiplication.  Bergeron and Li \cite{BL} prove
\begin{theorem}
If a graded algebra $A=\bigoplus_{n\geq 0}A_{n}$ over $\mathbb{C}$
satisfies Conditions (1)-(5) then $G(A)$ and $K(A)$ are graded dual
Hopf algebras.
\end{theorem}
In particular Theorem \ref{dimthm} applies to graded algebras which
satisfy Conditions (1)-(5).  Note that the dual Hopf algebras $G(A)$
and $K(A)$ come with distinguished bases consisting of the
isomorphism classes of simple and indecomposable projective modules.

\medskip
%%%%%%%%%%%%%%%%%%%%%%%%%%%%%%%%%%%
\section{From combinatorial Hopf algebras to dual graded graphs}
\label{sec:HNLS} This section recounts work of Fomin \cite{F},
  Hivert and Nzeutchap \cite{HN}, and Lam and Shimozono.  A {\it
graded graph} $\Gamma = (V,E,h,m)$ consists of a set of vertices
$V$, a set of (directed) edges $E \subset V \times V$, a height
function $h: V \to \{0,1,\ldots\}$ and an edge multiplicity function
$m: V \times V \to \{0,1,\ldots\}$. If $(v,u) \in E$ is an edge then
we must have $h(u) = h(v) + 1$. The multiplicity function determines
the edge set: $(v,u) \in E$ if and only if $m(v,u) \neq 0$.  We
assume always that there is a single vertex $v_0$ of height 0.

Let $\Z V = \bigoplus_{v \in V} \Z \cdot v$ be the free $\Z$-module
generated by the vertex set.  Given a graded graphs $\Gamma = (V, E,
h, m)$ we define up and down operators $U,D \colon\Z V \to \Z V$ by
$$
U_\Gamma(v) = \sum_{u \in V}m(v,u)\, u \ \ \ D_{\Gamma}(v) = \sum_{u
\in V}m(u,v)\, u
$$
and extending by linearity over $\Z$.  We will assume that $\Gamma$
is locally-finite, so that these operators are well defined.  A pair
$(\Gamma, \Gamma')$ of graded graphs with the same vertex set $V$
and height function $h$ is called {\it dual} with {\it differential
coefficient} $r$ if we have
$$ D_{\Gamma'} U_{\Gamma} - U_\Gamma D_{\Gamma'} =
r\,{\rm Id}.
$$
We shall need the following result of Fomin.  For a graded graph
$\Gamma$, let $f^v_{\Gamma}$ denote the number of paths from $v_0$
to $v$, where for two vertices $w, u \in V$, we think that there are
$m(w,u)$ edges connecting $w$ to $u$.
\begin{theorem}[Fomin~\cite{F}] \label{thm:F} Let $(\Gamma,\Gamma')$ be a pair of dual graded graphs
with differential coefficient $r$. Then
$$r^n n!=\sum_{v \colon h(v)=n}f^v_{\Gamma}f^v_{\Gamma'}.$$

\end{theorem}

Let $H_\bullet = \bigoplus_{n \geq 0} H_n$ and $H^\bullet =
\bigoplus_{n \geq 0} H^n$ be graded dual Hopf algebras over $\Z$
with respect to the pairing $\scal{\,.\,,.\,}: H_\bullet \times H^\bullet
\to \Z$. We assume that we are given dual sets of homogeneous free
$\Z$-module generators $\{p_\lambda \in H_\bullet\}_{\lambda \in
\Lambda}$ and $\{s_\lambda \in H^\bullet\}_{\lambda \in \Lambda}$,
such that all structure constants are non-negative integers.  We
also assume that $\dim(H_i) = \dim(H^i) < \infty$ for each $i \geq
0$ and $\dim(H_0) = \dim(H^0) = 1$, so that $H_0$ and $H^0$ are
spanned by distinguished elements the unit $1$.  Let us suppose we
are given non-zero homogeneous elements $\alpha \in H_1$ and $\beta
\in H^1$ of degree 1.

We now define a graded graph $\Gamma(\beta) = (V, E, h, m)$ where
$V=\{s_\lambda\}_{\lambda \in \Lambda}$ and $h\colon V\to{\mathbb
Z}$ is defined by $h(s_\lambda)=\deg(s_\lambda)$.  The map $m\colon
V\times V\to{\mathbb Z}$ is defined by
$$
m(s_\lambda,s_\mu)= \scal{p_\mu,\beta s_\lambda} =
\scal{\Delta(p_\mu), \beta\otimes s_\lambda}
$$
and $E$ is determined by $m$.  The graph $\Gamma(\beta)$ is graded because of  the assumption that $\beta$ has degree 1.  Similarly,
we define a graded graph $\Gamma'(\alpha)=(V',E',h',m')$ where $V' =
V$, $h' = h$, and
$$
m'(s_\lambda,s_\mu)= \scal{\alpha \, p_\lambda, s_\mu} = \scal{
\alpha \ot  p_\lambda, \Delta(s_\mu)}.$$  The following theorem is
due independently to Hivert and Nzeutchap \cite{HN} and Lam and
Shimozono (unpublished).

\begin{theorem}\label{thm:HNLS}
The graded graphs $\Gamma = \Gamma(\beta)$ and $\Gamma' =
\Gamma'(\alpha)$ form a pair of dual graded graphs with {\it
differential coefficient} $\scal{\alpha,\beta}$.
\end{theorem}
\begin{proof}
We identify $\Z V$ with $H^\bullet$ and note that $U_\Gamma(x) =
\beta \, x$ where $x \in H^\bullet$ and we use the multiplication in
$H^\bullet$. Also,
\begin{align*}
D_{\Gamma'}(x) &= \sum_{\mu \in \Lambda} \scal{\alpha \ot p_\mu,
\Delta x} \, s_\mu = \sum \scal{\alpha,x^{(1)}}\, x^{(2)}.
\end{align*}
where $\Delta x = \sum x^{(1)} \otimes x^{(2)}$. Now observe that by
our hypotheses on the degree of $\alpha$ and $\beta$ they are
primitive elements: $\Delta\alpha = 1 \otimes \alpha + \alpha
\otimes 1$ and $\Delta\beta = 1 \otimes \beta + \beta \otimes 1$. We
first calculate
\begin{align*}
\scal{\alpha,\beta \, x} =\scal{\Delta\alpha,\beta \otimes x}
=\scal{1,\beta}\scal{\alpha,x}+\scal{\alpha,\beta}\scal{1,x}
=\scal{\alpha,\beta}\scal{1,x}
\end{align*}
and then compute
\begin{align*}
D_{\Gamma'} U_{\Gamma}(x) & = D_{\Gamma'}(\beta x)\\
& = \sum \left(\scal{\alpha,\beta \, x^{(1)}}\,x^{(2)} +
\scal{\alpha,x^{(1)}}\,\beta\,x^{(2)}\right) \\
&= \scal{\alpha,\beta}x +U_{\Gamma}D_{\Gamma'}(x)
\end{align*}
where to obtain $\scal{\alpha,\beta}x$ in the last line we use
$\Delta x = 1 \ot x + \text{terms of other degrees}.$
\end{proof}

%%%%%%%%%%%%%%%%%%%%%%%%%%%%%%%%%%%
\section{Proof of Theorem \ref{dimthm}}
We are given a graded algebra $A=\bigoplus_{n\geq 0}A_{n}$ over
$\mathbb{C}$ with multiplication $\rho$ satisfying Conditions (1),
(2) and (3).  Moreover we assume that the two Grothendieck groups
$G(A)$ and $K(A)$ form a pair of graded dual Hopf algebras as in
Section \ref{sec:BL}.  Under these assumptions we show that 
$$ \dim(A_n)=r^n n! $$ where $r=\dim(A_1)$.

Let $H^{\bullet}=G(A)$ and $H_{\bullet}=K(A)$.  Let $\{s^{(1)}_1 =
[S^{(1)}_1],\ldots,s^{(1)}_t = [S^{(1)}_t]\}$ and $\{p^{(1)}_1=
[P^{(1)}_1],\ldots,p^{(1)}_t = [P^{(1)}_t]\}$ denote the isomorphism
classes of simple and indecomposable projective $A_1$-modules, so
that $H^1=\bigoplus_{i=1}^{t}\mathbb{Z}s^{(1)}_i$ and
$H_1=\bigoplus_{i=1}^{t}\mathbb{Z}p^{(1)}_i$.  Define
$a_i=\dim(S^{(1)}_i)$ and $b_i=\dim(P^{(1)}_i)$ for $1\leq i\leq t$.
We set for the remainder of this paper
    $$\alpha=\sum_{i=1}^ta_ip^{(1)}_i\in H_{1}\qquad \hbox{and}\qquad \beta=\sum_{i=1}^tb_is^{(1)}_i\in H^{1}.$$

Since $A_0\cong \mathbb{C}$, we let $s^{(0)}_1$ (respectively,
$p^{(0)}_1$) be the unique simple (respectively, indecomposable
projective) module representative in $H^0$ (respectively, $ H_0$).
Similarly, let $\{s^{(n)}_i = [S^{(n)}_i]\}$ be all isomorphic
classes of simple $A_n$-modules and $\{p^{(n)}_i = [P^{(n)}_i]\}$ be
all isomorphism classes of indecomposable projective $A_n$-modules.
The sets $\bigcup_{n\ge 0}\{s^{(n)}_i\}$ and $\bigcup_{n\ge
0}\{p^{(n)}_i\}$ form dual free $\Z$-module bases of $H^{\bullet}$
and $H_{\bullet}$.

Now define $\Gamma = \Gamma(\beta)$ and $\Gamma' = \Gamma'(\alpha)$
as in Section \ref{sec:HNLS}.

\begin{lemma}\label{lem:dimpath}
We have
$$
f^{s^{(n)}_j}_{\Gamma} = \dim P^{(n)}_j \ \ \text{and} \ \
f^{s^{(n)}_j}_{\Gamma'} = \dim S^{(n)}_j.
$$
\end{lemma}
\begin{proof}
We have
$$
m(s^{(n-1)}_i,s^{(n)}_j)=\sum_{l=1}^t b_lc_l,
$$
where $c_l$ is the number of copies of the indecomposable projective
module $P^{(1)}_l\otimes P^{(n-1)}_i$ as a summand in
Res$_{A_{1}\otimes A_{n-1}}^{A_n}P^{(n)}_j$.  Note that $s^{(0)}_1$
is the unit of $H^{\bullet}$ and $m(s^{(0)}_1,s^{(1)}_i)=b_i=\dim
P^{(1)}_i$ for all $1\leq i\leq t$. The dimension of an
indecomposable projective module $P^{(n)}_j$ is given by
  $$\dim P^{(n)}_j = \sum_{i,l} c_l \dim  \left(P^{(1)}_l\otimes P^{(n-1)}_i\right)= \sum_i m(s^{(n-1)}_i,s^{(n)}_j)\dim P^{(n-1)}_i.
  $$
By induction on $n$,  we deduce that $\dim P^{(n)}_j$ is the number
of paths from $s^{(0)}_1$ to $s^{(n)}_j$ in $\Gamma$.  The claim for
$\Gamma'$ is similar.
\end{proof}

For any  finite dimensional  algebra $B$ let $\{S_{\lambda}\}_{\lambda}$ be a complete set of simple
 $B$-modules. For each $\lambda$ let $P_\lambda$ be the projective cover of $S_\lambda$. It is well known (see~\cite{CR}) that we can find minimal idempotents $\{e_i\}$ such that $B=\bigoplus Be_i$ where each $Be_i$ is isomorhpic to a $P_\lambda$. Moreover, the quotient of $B$ by its radical shows that the multiplicity of $P_\lambda$ in $B$ is equal to $\dim S_\lambda$. This implies the following lemma.

\begin{lemma}\label{dimreplem} Let $B$ be a finite dimensional  algebra and $\{S_{\lambda}\}_{\lambda}$ be a complete set of simple $B$-modules.
$$\dim B=\sum_{\lambda} (\dim P_{\lambda})(\dim S_{\lambda}),$$
where $P_{\lambda}$ is the  projective cover of $S_{\lambda}$.
\end{lemma}

By Lemma \ref{dimreplem}, we have $r= \sum_{i=1}^t a_ib_i=
\scal{\alpha,\beta}$.  By Theorem \ref{thm:HNLS} we may apply
Theorem \ref{thm:F} to $(\Gamma,\Gamma')$.  Using Lemma
\ref{dimreplem} and Lemma \ref{lem:dimpath}, Theorem \ref{thm:F}
says
$$
\dim(A_n) = \sum_i (\dim P^{(n)}_{i})(\dim S^{(n)}_{i}) = \sum_i
f_\Gamma^{s_i^{(n)}}\,f_{\Gamma'}^{s_i^{(n)}} = r^n n!.
$$

\begin{remark} If the tower consists of semisimple algebras $A_i$
then $\Gamma = \Gamma'$ so we obtain a {\it self-dual graph}.  In
this case the graph would be a weighted version of a {\it
differential poset} in the sense of Stanley \cite{S}.  If
furthermore the branching of irreducible modules from $A_n$ to
$A_{1} \otimes A_{n-1}$ is multiplicity free then we get a true
differential poset.
\end{remark}

\begin{remark} The Hopf algebras $H^\bullet$ and $H_\bullet$ are not
in general commutative and co-commutative.  Thus in the definitions
of Section \ref{sec:HNLS} we could have obtained a different pair of
dual graded graphs by setting $m(s_\lambda,s_\mu)= \scal{p_\mu,
s_\lambda\, \beta}$  or $m'(s_\lambda,s_\mu)= \scal{p_\lambda \,
\alpha, s_\mu}$.
\end{remark}

%%%%%%%%%%%%%%%%%%%%%%%%%%%%%%%%%%%
%%%%%%%%%%%%%%%%%%%%%%%%%%%%%%%%%%%
%\section{Concluding remarks}
%
%[Do we say anything more?]
%%%%%%%%%%%%%%%%%%%%%%%%%%%%%%%%%%%
%%%%%%%%%%%%%%%%%%%%%%%%%%%%%%%%%%%
%%%%%%%%%%%%%%%%%%%%%%%%%%%%%%%%%%%

\end{document}